\documentclass[5p]{elsarticle}
\usepackage{cite}
\usepackage{amsmath,amssymb,amsfonts}
\usepackage{algorithmic}
\usepackage{graphicx}
\usepackage{textcomp}
\usepackage{booktabs}
\usepackage{url}
\usepackage[skip=4pt]{caption}

\begin{document}
\begin{frontmatter}

\title{Hyperellipsoid Density Sampling:\\Exploitative Sequences to Accelerate High-Dimensional Numerical Optimization}

\author{Julian G. Soltes}

\affiliation{organization={Regis University},
city={Denver},
state={Colorado},
country={USA}}

\begin{abstract}
The curse of dimensionality remains a persistent challenge in modern optimization problems. Expanding the search space into higher dimensions exponentiates the sparsity of uniform sample distributions, rendering traditional quasi-Monte Carlo (QMC) sequences increasingly inefficient. This paper introduces a non-uniform sampling strategy to accelerate high-dimensional optimization.

This method, Hyperellipsoid Density Sampling (HDS), generates samples as hyperellipsoids overlapping throughout the parameter space. Utilizing a series of unsupervised learning techniques, a non-uniform sequence is generated to exploit the interior regions of the hypervolume. If prior information about optima is known, the distribution can be biased towards the known regions, making HDS versatile for many numerical applications.

HDS was evaluated against Sobol, a highly uniform QMC sampling method, using differential evolution (DE) on the challenging and widely benchmarked set of 29 CEC2017 test functions. The results show statistically significant improvements in final solution geometric mean error ($p<0.05$), with average performance gains ranging from $37\%$ in $10D$ to $11\%$ in $100D$. This paper demonstrates the efficacy of HDS as an exploitative alternative to uniform QMC sampling.
\end{abstract}


\end{frontmatter}

\section{Introduction}
Traditional random or quasi-Monte Carlo uniform sampling methods aim to maximize coverage of the search space with minimal discrepancy. In optimization problems, where the global optimum often resides in interior regions of the parameter space, these methods can become inefficient due to increasing sparsity \citep{Chen2025}. Specifically in high dimensions, sampling uniformly across the entire hypercube consumes computational resources exploring the typically irrelevant corners of the space, where the majority of the hypervolume is concentrated.

\begin{figure}[b!]
\centering
\includegraphics[width=0.4\textwidth]{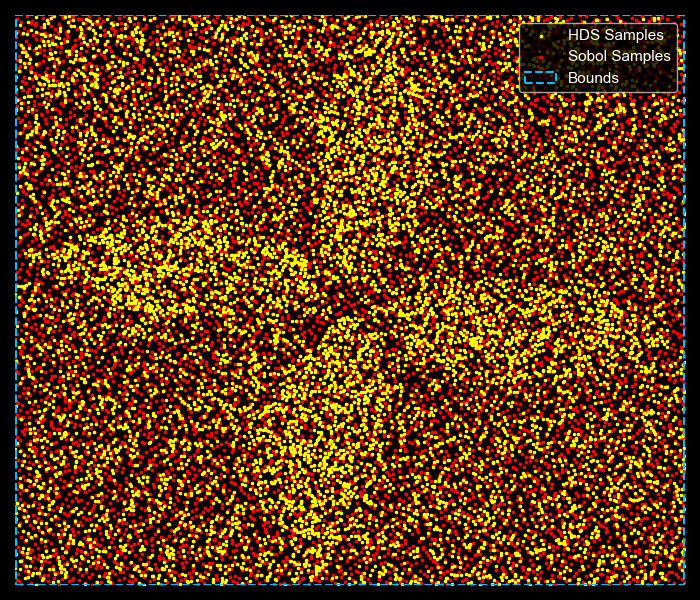}
\caption{HDS (yellow) and comparison Sobol sequence (red) for $N=10,000$ in $2D$.}
\label{fig:hds_vs_sobol_sequences}
\end{figure}

Various approaches have been explored to address this scaling inefficiency, such as modifying traditional optimizers and employing data preprocessing methods \citep{Hinz2018}.
This paper presents an alternative strategy: ellipsoidal sample sequences (Fig. \ref{fig:hds_vs_sobol_sequences}) that densely populate inner regions of high-dimensional spaces.

\subsection{Hyperellipsoid Density Sampling}
The proposed strategy, Hyperellipsoid Density Sampling (HDS), generates its distributions by focusing a large initial low-discrepancy sequence down to the target sample size. By densely populating interior regions of the parameter space, HDS accounts for the geometric sparsity that renders uniform sampling inefficient in high dimensions.

This is achieved by defining multiple hyperellipsoids over the parameter space, derived from a series of unsupervised learning methods performed on an initial Sobol QMC sequence. The result is a non-uniform sample sequence concentrated in areas more likely to contain the global optimum compared to sparse uniform sampling.

\section{Methodology}
\subsection{Experimental Methods}
This study evaluates the performance of different initial sampling strategies using differential evolution (DE) as the test
optimizer. DE is an effective and widely used algorithm for non-differentiable, non-linear, multi-modal, and high-dimensional problems \citep{Stokes2020}. Each trial is conducted using the \texttt{Scipy} Differential Evolution implementation, with default hyperparameters ($\text{strategy}=\text{'best1bin'}$, $F=(0.5,1)$, $CR=0.7$) for baseline use-case performance.

\subsubsection{Optimization Trials}
The experimental trials were run against the 29 CEC2017 benchmark test functions, each possessing distinct characteristics to challenge optimizers and are scalable to 10, 30, 50, and 100 dimensions ($D \in  \{10,30,50,100\}$) \citep{Naser2025}. As such, these are the four dimensionalities tested, with bounds of [-100,100]. 

50 trials were conducted for two sets of sample sizes, as well as for each combination of test function and dimension, totaling 11,400 DE runs. The sample sizes tested are 64 and 1000; the lower value is chosen for sparser coverage with a power of 2 for Sobol uniformity \citep{Joe2008}, the higher value represents a more typical optimization problem. The number of evolutionary iterations is fixed at $n_\text{iter} = 100$, corresponding to $F_\text{eval} = N \times 100 \implies F_\text{eval} \in \{6400, 100000\}$ total function evaluations in each trial. This relatively low $F_\text{eval}$ highlights the practical benefits of initializing a computationally limited search using HDS compared to uniform methods.

\subsubsection{Comparison Metrics}
The comparison metric used for the analysis is the final solution fitness errors for HDS compared to Sobol. Sobol was chosen as the comparison metric due to its high uniformity and low discrepancy in high dimensions \citep{Atanassov2022}, as well as its use within the HDS algorithm. Latin Hypercube was considered as an alternative, but does not ensure the same low discrepancy as Sobol in high dimensions \citep{Wang2007}. To ensure a fair comparison, the geometric mean is used to calculate the average improvement factors for each dimension across all test functions. 

The optimization run times and sample generation times between the two sampling methods are compared similarly, to assess the separate processing times.

The discrepancies of the resulting HDS sequences are additionally compared against Sobol using L2-Star and Centered L2 metrics.

\subsection{Sample Sequences}
The HDS method generates $N$ total samples $\mathbf{x} \in \mathbb{R}^D$ within a $D$-dimensional search space (defined by the bounds \textbf{B}), using a geometric approach to focus sample density. This relies purely on the structure of an initial Sobol sequence, requiring zero objective function evaluations. 

\begin{figure}[t!]
\centering
\includegraphics[width=0.5\textwidth]{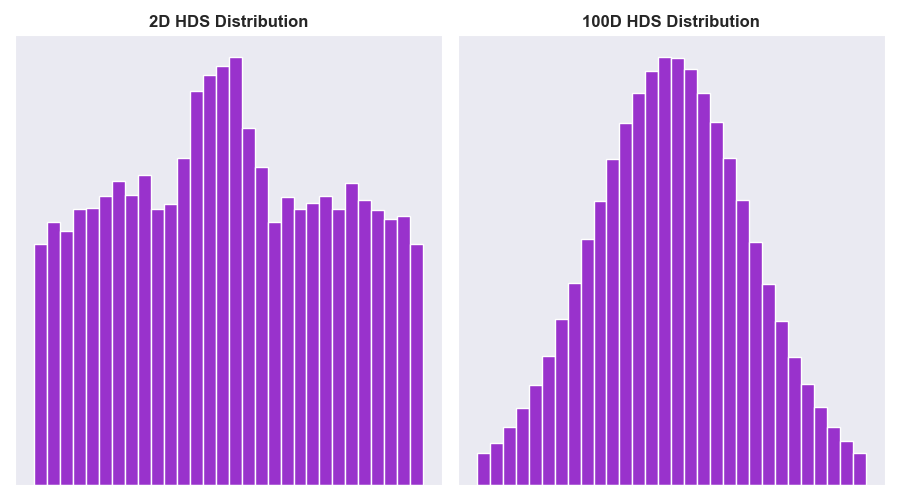}
\caption{Sample distributions for HDS sequences in $2D$ (left) and $100D$ (right).}
\label{fig:hds_distributions}
\end{figure}

The resulting sample distribution is highly non-uniform, increasingly resembling a normally projected distribution as dimensions scale higher (Fig. \ref{fig:hds_distributions}). This accounts for the curse of dimensionality by replacing uniform coverage with an exploitative density; if prior information is known, the distribution may be additionally biased by applying Gaussian weights.

\subsubsection{Define Normalized Space}
The D-dimensional search bounds, $\mathbf{B}$, are normalized to a unit hypercube $[0,1]^D$ to simplify clustering and Principal Component Analysis (PCA).

\subsubsection{Generate Initial QMC Sequence}
A large initial set of $N_{\text{init}}$ scrambled Sobol samples is generated, serving as a maximally uniform and minimally discrepant structural distribution of the parameter space. $N_{\text{init}}$ is set to the smallest power of 2 greater than or equal to $200 \cdot D$, ensuring sufficient coverage that scales efficiently with dimensionality as $N_{\text{init}} = 2^{\lceil \log_2 (200 D) \rceil}$. This value is capped at $2^{15}$ by default for efficiency.

These values balance computational efficiency with spatial discrepancy for high sample sizes and dimensions, as well as ensuring power-of-2 sample sizes to minimize Sobol sampling artifacts. This initial set serves as the basis for determining the geometry of the search space.

\subsubsection{Apply Gaussian Weights (Optional)}
\begin{figure}[t!]
\centering
\includegraphics[width=0.4\textwidth]{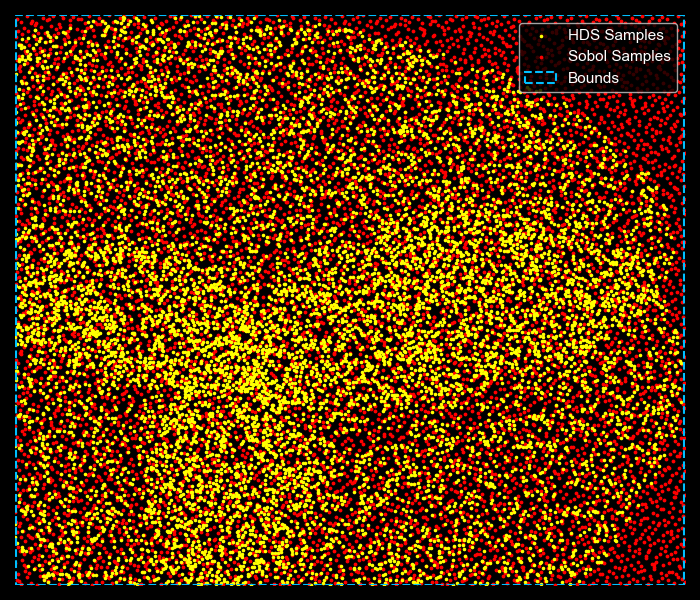}
\caption{$2D$ sequence with weights applied towards $(0.25, 0.25)$, for bounds $\mathbf{[0,1]^D}$.}
\label{fig:hds_weighted}
\end{figure}

If prior information (such as target distribution or optimal solution region) is known, Gaussian weights may be applied to the initial QMC samples to influence their distributions. As a result, the final sample locations are shifted towards the desired region. An example of an HDS sample sequence with weights applied is shown in Fig. \ref{fig:hds_weighted}.

\subsubsection{Initial Centroid Estimation}
The initial $N_{\text{init}}$ samples are clustered using MiniBatchKMeans to reduce the data into a smaller, computationally efficient set of $K_{\text{init}}$ centroids (defaulted at 100), simplifying the subsequent hierarchical analysis. MiniBatchKMeans is used in place of standard KMeans to minimize sample generation time.

The value $K_{\text{init}}=100$ was chosen as a stable value, as it provides the optimal balance between spatial information and computational efficiency. This value consistently generates clear KMeans and Hierarchical clustering results to sufficiently cover the parameter space while minimizing the high overhead of the linkage calculation. As dimensionality approaches $1000D$, this default value is scaled down proportionally, ensuring the final sample distribution remains minimally affected.

\subsubsection{Identify Number of Ellipsoids}
\begin{figure}[t!]
\centering
\includegraphics[width=0.4\textwidth]{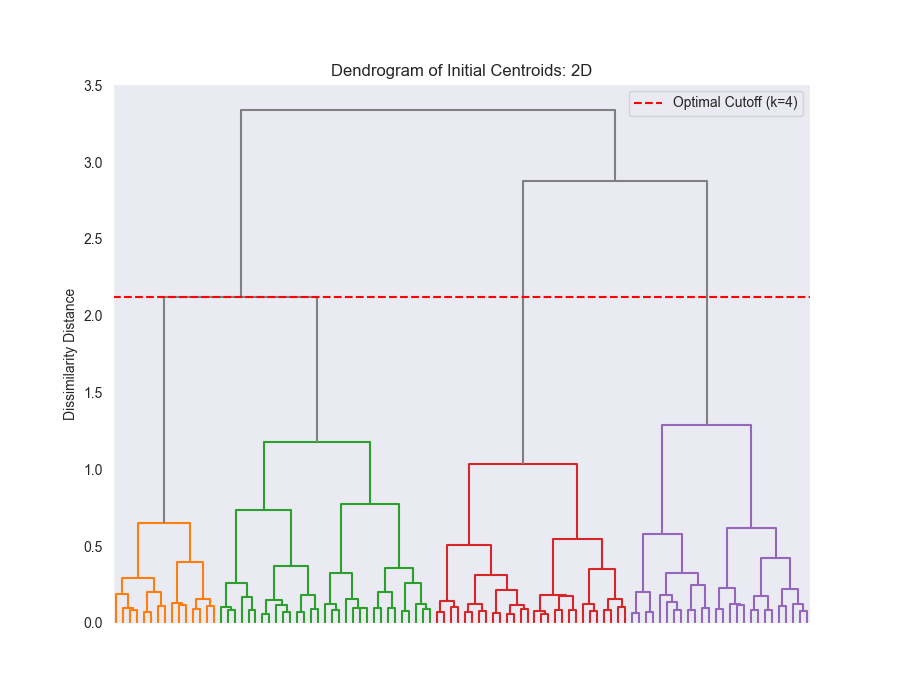}
\caption{Dendrogram used to identify number of ellipsoids for a representative $2D$ sample sequence.}
\label{fig:dendrograms}
\end{figure}

Agglomerative Hierarchical Clustering (AHC) is performed on the $K_{\text{init}}$ centroids. The optimal number of clusters, $K$, is determined by analyzing the dendrogram's linkage matrix (Fig. \ref{fig:dendrograms}) and selecting a cut-off distance $\delta$ corresponding to the largest dissimilarity between the clusters. This distance threshold determines the final number of clusters, which in turn defines the number of ellipsoids used in the subsequent calculations.

Due to the stochastic nature of the initialization process, the optimal cluster count $K$ may vary across multiple sample generations with identical input parameters. The number of ellipsoids can also be provided as an input parameter, and the AHC calculations will be skipped.

\subsubsection{Identify Ellipsoid Centroids}
A final MiniBatch-KMeans fit, using the optimal cluster count $K$, assigns each of the $N_{\text{init}}$ samples to one of $K$ clusters, yielding the $\mathbf{c}_k$ cluster centers (ellipsoid origins) and the sample counts $n_k$ for clusters $k \in \{1,2, ..., K \}$.

\subsubsection{Derive Ellipsoid Geometry}
The shape and orientation of the $K$ hyperellipsoids are defined using PCA on the samples within each cluster.

For each cluster $k$, the origin of the hyperellipsoid is set to the centroid $\mathbf{c}_k$. PCA is performed on the associated samples, $\mathbf{s}_k$. The principal components define the orientation of the ellipsoid axes, $\mathbf{P}$. The length of the $i$-th semi-axis, $\sigma_{i,k}$, is determined by the calculated variance $\text{Var}_{i,k}$ along that component:

$$\sigma_{i,k} = \sqrt{\text{Var}_{i,k} + \epsilon}$$

where $\epsilon$ is a small constant added for numerical stability in the case of zero variance.

\subsubsection{Allocate Samples}
The total number of required samples, $N$, is distributed among the $K$ hyperellipsoids based on the relative density of their initial clusters. The number of samples allocated to cluster $k$, $N_k$, is
$$N_k = \text{round}\left(N \frac{n_k}{\sum_{i=1}^{K} n_i}\right)$$

where $n_k$ is the number of initial QMC samples belonging to cluster $k$.

\subsubsection{Initialize Unit Hyperspheres}
For each cluster $k$, a set of direction vectors $\mathbf{u}$ on the surface of the unit sphere is generated using the Marsaglia polar method. Each sample vector $\mathbf{s}$ is normalized to a unit length $\mathbf{\hat{u}}$ by dividing by its L2-norm:

$$\mathbf{\hat{u}} = \frac{\mathbf{s}}{\|\mathbf{s}\|}$$

where $\mathbf{s} \sim \mathcal{N}(\mathbf{0}, \mathbf{I}_D)$ is an initial vector with $D$ independent standard normal components.

\subsubsection{Radial Scaling Factor}
A global scaling factor $\lambda$ is applied to all hyperellipsoid axes to control the spread of the sampling points. This factor is calculated by scaling the square root of the $\chi^2$ critical value with an empirically tuned, dimension-variant constant $C_D$. The confidence level of the $\chi^2$ critical value is chosen as $\alpha=0.9999$ to capture the majority of each cluster's variance. The empirically tuned constant $C_D$ is a function of the dimension $D$:

$$C_D = 0.55 - 0.01 \log_2(D)$$
The final (default) global scaling factor $\lambda$ is then defined as:
$$\lambda = C_D \sqrt{\chi^2_{\alpha, D}}$$

\begin{figure}[t!]
\centering
\includegraphics[width=0.4\textwidth]{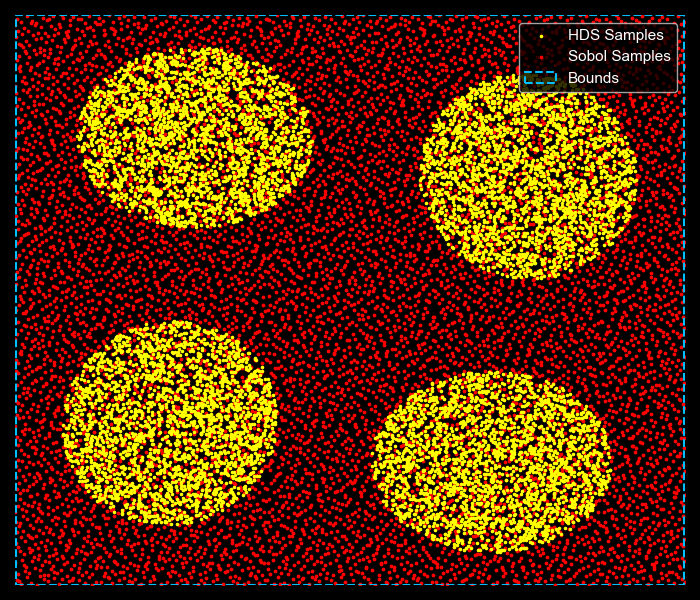}
\caption{HDS samples in $2D$ with $K=4$ ellipsoids, using a reduced radial scaling factor.}
\label{fig:4ellipsoids}
\end{figure}

This scaling is critical; a factor that is too small leads to poor coverage, and too large leads to a sample high rejection rate (described below). A $2D$ example is shown in Fig. \ref{fig:4ellipsoids} generated with a reduced scaling factor $\lambda$ to easily visualize the geometry.

\subsubsection{Radial Scaling}
To ensure a uniform distribution within each ellipsoid, a $1D$ Sobol sequence $q \in [0,1]$ is generated and projected along the unit direction vector $\mathbf{\hat{u}}$. Taking $q^{1/D}$ ensures that the radial distances are uniformly distributed throughout the hypervolume, scaled by the global scaling factor $\lambda$.

$$\mathbf{u} = \left(q^{1/D} \lambda \right) \mathbf{\hat{u}}$$

\subsubsection{Hyperellipsoid Transformation}
The scaled unit sphere samples, now a matrix $\mathbf{S}$, are transformed into the final hyperellipsoid by applying the variance-based axis lengths and rotating them back to the original parameter space. The full transformation, including scaling and rotation, is defined by matrix multiplication:

$$\mathbf{X} = \mathbf{S} \mathbf{\Sigma} \mathbf{P}^{\top} + \mathbf{C}_k$$

where:

\begin{itemize}
    \item $\mathbf{X}$ is the matrix of ellipsoid samples.
    \item $\mathbf{S}$ is the matrix of scaled unit direction vectors.
    \item $\mathbf{\Sigma}$ is the diagonal matrix of scaled axis lengths.
    \item $\mathbf{P}^{\top}$ is the PCA rotation matrix.
    \item $\mathbf{C}_k$ is the vector of cluster centroids $\mathbf{c}_k$, added to every row of the scaled and rotated matrix.
\end{itemize}

\subsubsection{Void Filling}
Samples generated outside the normalized search bounds $[0,1]^D$ are rejected. A large number of initial samples is generated to compensate for this expected rejection rate.

If the count of valid samples is less than the target $N$, the remaining points ($N_{\text{void}}$) are generated using an adaptive void-filling strategy. This method identifies the most sparse regions using a BallTree structure to optimize a nearest-neighbor search. The resampled points are then generated using a truncated normal distribution, centered on existing samples in these sparse regions. This ensures the total sample count is met while providing explorative coverage to areas not fully contained by the hyperellipsoids, which was found to significantly improve performance over non-adaptive void-filling methods.

\subsubsection{Denormalization}
The resulting samples $\hat{\mathbf{H}}$, currently normalized between $[0,1]$, are scaled element-wise back to the original search bounds $\mathbf{B}$ using the minimum $\mathbf{B}_{\min}$ and span $\Delta \mathbf{B}$. This returns the final HDS sequence $\mathbf{H}$.

\begin{equation}
\mathbf{H} = \hat{\mathbf{H}} \odot \Delta \mathbf{B} + \mathbf{B}_{\min}
\end{equation}

\section{Results}
\subsection{Solution Fitness}
The final solution fitness being minimized is used as the metric to compare optimization efficacy of HDS compared to QMC sampling. The analysis is performed on the differential evolution results from the 50 independent trials, for each of the 2 sample sizes, 4 dimensions and 29 CEC2017 benchmark test functions. The experimental trials on the benchmark functions demonstrated that HDS consistently outperformed the Sobol QMC in all tested dimensions. The 95\% confidence intervals show similar performance between the sampling methods in lower dimensions ($10D$, $30D$), whereas in higher dimensions ($50D$, $100D$) the $CI_{95}$ ranges do not fall below 1.0x improvement factors. This highlights the consistent improvements of HDS in higher-dimensional problems.

The $N=1000$ trials saw an overall geometric mean improvement factor of 1.15x, implying that HDS averages 15\% improved solution quality across all dimensions tested. The $N=64$ trials showed similar but more nuanced improvements, lacking statistical significance for $10D$ and $100D$; the $10D$ case specifically showed a 3\% decrease in solution quality compared to Sobol, with highly variable performance leading to a $p$-value of $0.15$.

As shown in Table \ref{tab:optimization_results}, HDS achieved improvements in geometric mean solution quality aggregated across all test functions. The initial geometric bias of HDS accelerates convergence during the critical early phase of the search, shown in by the anytime performance \citep{hansen2022anytime} in Fig. \ref{fig:anytime_performance} in the appendix. Statistically significant improvements were measured across all tested dimensions in the $N=1000$ trials, as well as most dimensions (excluding $10D$) in the $N=64$ trials.

\begin{table}[h]
\scriptsize
\centering
  \caption{Final Solution Quality (Minimizing)}
\label{tab:optimization_results}
  \setlength{\tabcolsep}{4pt}
  \begin{tabular}{r r r c r r r}
    \toprule
    $N$ & $D$ & $R_\text{err}$ & $p$-Value & $\text{CI}_{95}$ & HDS Err. & Sobol Err. \\
    \midrule
    \midrule
    $1000$ & $10$ & $1.37$ & $3 \times 10^{-2}$ & (0.92, 2.05) & $6.0 \times 10^{-1}$ & $8.3 \times 10^{-1}$ \\
    $1000$ & $30$ & $1.03$ & $1 \times 10^{-4}$ & (0.89, 1.19) & $6.7 \times 10^{3}$ & $6.9 \times 10^{3}$ \\
    $1000$ & $50$ & $1.11$ & $4 \times 10^{-5}$ & (1.06, 1.15) & $5.2 \times 10^{4}$ & $5.7 \times 10^{4}$ \\
    $1000$ & $100$ & $1.11$ & $4 \times 10^{-6}$ & (1.06, 1.15) & $5.8 \times 10^{5}$ & $6.4 \times 10^{5}$ \\
    \midrule
    $64$ & $10$ & $0.97$ & $0.15$ & (0.67, 1.43) & $8.9 \times 10^{0}$ & $8.7 \times 10^{0}$ \\
    $64$ & $30$ & $1.17$ & $3 \times 10^{-5}$ & (1.05, 1.30) & $1.1 \times 10^{4}$ & $1.2 \times 10^{4}$ \\
    $64$ & $50$ & $1.12$ & $1 \times 10^{-6}$ & (1.07, 1.18) & $7.5 \times 10^{4}$ & $8.4 \times 10^{4}$ \\
    $64$ & $100$ & $1.08$ & $0.16$ & (1.04, 1.13) & $8.5 \times 10^{5}$ & $9.2 \times 10^{5}$ \\
    \bottomrule
  \end{tabular}
\end{table}

\subsection{Computational Overhead}
The multiple machine learning fits and backend QMC generations naturally lead to slower sample generation than a single deterministic QMC sequence. Generation time $t_g$ increases linearly with sample size and dimension (Table. \ref{tab:generation_time}). The optimization run time ratios $R_\text{time}$ between HDS and Sobol are highly influenced by $t_g$ in low $F_\text{eval}$ domains, highlighted in Table \ref{tab:runtime_results}. The arithmetic average run time ratio is found to be $\bar{R}_\text{time}=0.951$ ($p < 1 \times 10^{-291}$), indicating $\sim 5\%$ total increase in execution time.

Due to the scaling of sample generation time, HDS sampling in extremely high dimensions ($D>1000$) and sample sizes ($N>10,000$) may become computationally prohibitive. However, in these domains, the increased overhead of initializing a single sequence may be small compared to the total optimization time.

\begin{table}[h]
    \centering
    \caption{Sample Generation Time}
    \label{tab:generation_time}
    \begin{tabular}{r r r c}
      \toprule
      $N$ & $D$ & Avg. (s) & $\text{CI}_{95}$ \\
      \midrule
      \midrule
      10 & 100 & 1.59 & (1.46, 1.72) \\
      100 & 100 & 1.40 & (1.33, 1.47) \\
      1,000 & 100 & 1.46 & (1.40, 1.51) \\
      10,000 & 100 & 1.70 & (1.62, 1.78) \\
      100,000 & 100 & 2.54 & (2.46, 2.63) \\
      \midrule
      1,000 & 1 & 0.17 & (0.16, 0.18) \\
      1,000 & 10 & 0.36 & (0.32, 0.39) \\
      1,000 & 50 & 1.12 & (1.06, 1.18) \\
      1,000 & 100 & 1.67 & (1.58, 1.75) \\
      1,000 & 1,000 & 6.83 & (6.66, 7.00) \\
      \bottomrule
    \end{tabular}
\end{table}

\begin{table}[h]
\centering
\caption{Run Time Ratios}
    \label{tab:runtime_results}
\begin{tabular}{r r c r}
      \toprule
      $N$ & $D$ & $R_\text{time}$ & $p$-Value \\
      \midrule
      \midrule
      $1000$ & $10$ & $0.96$ & $4 \times 10^{-12}$ \\
      $1000$ & $30$ & $0.96$ & $1 \times 10^{-51}$ \\
      $1000$ & $50$ & $0.96$ & $7 \times 10^{-87}$ \\
      $1000$ & $100$ & $0.93$ & $1 \times 10^{-170}$ \\
      \midrule
      $64$ & $10$ & $0.64$ & $6 \times 10^{-195}$ \\
      $64$ & $30$ & $0.69$ & $1 \times 10^{-186}$ \\
      $64$ & $50$ & $0.65$ & $9 \times 10^{-215}$ \\
      $64$ & $100$ & $0.49$ & $6 \times 10^{-233}$ \\
      \bottomrule
    \end{tabular}
\\
\footnotesize{$^*R_\text{time} < 1$ indicates slower HDS performance.}
\end{table}

\subsection{Discrepancies}
HDS inherently sacrifices uniformity due to its exploitative nature. This is explored by quantifying the discrepancies of the HDS and pure Sobol sample sequences. 

The HDS L2-Star discrepancies are significantly higher than those of Sobol; Centered L2 values, however, are consistently lower than Sobol due to the tighter clustering of samples. Together, this implies high local sample density with lower global uniformity across the full parameter space. For a set of 30 representative $100D$ trials, the discrepancy comparisons using L2-Star and Centered L2 metrics are shown in Table \ref{tab:discrepancies}.

\begin{table}[h]
  \centering
  \caption{Discrepancies}
  \label{tab:discrepancies}
  \setlength{\tabcolsep}{6pt} 
  \begin{tabular}{r c c c c}
    \toprule
     & \multicolumn{2}{c}{\textbf{L2-Star ($\times 10^{-17}$)}} & \multicolumn{2}{c}{\textbf{Centered-L2 ($\times 10^4$)}} \\
    \cmidrule(lr){2-3} \cmidrule(lr){4-5}
    \textbf{$N$} & HDS & Sobol & HDS & Sobol \\
    \midrule
    \midrule
    $10$ & $48$ & $0.000023$ & $5500$ & $32000$ \\
    $100$ & $34$ & $0.000300$ & $180$ & $5500$ \\
    $500$ & $37$ & $0.000060$ & $19$ & $1000$ \\
    $1000$ & $26$ & $0.000170$ & $8$ & $490$ \\
    \bottomrule
  \end{tabular}
\end{table}

\section{Ease of Use \& Implementation}
HDS is designed to generate samples needing only the number of samples and parameter bounds, making it seamlessly integratable within existing workflows. It uses a small set of standard library imports for minimal external dependency.

HDS is implemented in PyPi as part of the \texttt{hdim-opt} Python package, following the conventions of Scipy’s QMC module. The open-source code containing the experimental notebooks can additionally be found at \url{https://github.com/jgsoltes/Hyperellipsoid-Density-Sampling} for full reproducibility.

\newpage
\section{Conclusion}
Hyperellipsoid Density Sampling (HDS) is introduced as an effective sampling method for high-dimensional optimization. These sequences densely populate the interior regions of the parameter space, as an exploitative alternative to uniform Quasi-Monte Carlo (QMC) sampling.

Biasing the distributions via Gaussian weights allows the inclusion of a priori information to further exploit regions of known interest. This extends the utility of HDS beyond traditional numerical optimization. 

The experimental results show that HDS leads to statistically significant improvements in final solution quality compared to standard uniform QMC sequences. Across all dimensions tested, the geometric mean improvement was $15\%$ compared to Sobol, with only a $5\%$ increase in total computation time. The efficacy of HDS highlights the value of using hyperellipsoidal distributions to accelerate an initial search, beneficial for applications such as data modeling, sensitivity analysis, and deep reinforcement training.

\section{Future Works}
Future works explore the performance of HDS against a wider range of optimization algorithms and test functions, as well as sample size-dimensionality permutations, to validate its robustness. 

\subsection{Hyperparameters}
Many hyperparameters within the HDS generation are derived empirically. A more rigorous mathematical framework may be explored to utilize the high-dimensional geometry for more efficient sample generation. Further work focuses on rigorously defining parameters such as initial QMC sample size, number of initial estimation clusters, as well as miscellaneous constants and scaling factors. This includes optimizing the algorithm’s sample generation time, especially for high sample counts and dimensions $D \gg 100$ and $N \gg 10,000$. Using covariance matrix diagonals rather than PCA demonstrates slight improvements, at the cost of degraded experimental solution quality.

\subsection{Weights}
The trials in this study were performed without the application of Gaussian weights. Evaluating the impact of these weights, specifically on problems where prior information about the optimal solution location is known, would validate the effectiveness of this weighting strategy. Hybrid sampling-optimization strategies can additionally be explored, including a first optimization initialized with a uniform distribution, followed by a second optimization initialized with HDS to exploit the previously identified region. This process can be repeated ad infinitum to converge sample sequences around the a priori regions of interest.

\section*{Acknowledgements}
The author would like to thank Dr. Ksenia Polson, Dr. Kellen Sorauf, Dr. Mike Busch, and Monet Morris for their invaluable feedback and support throughout this study. It would not have been possible without the insights gained from each of our discussions.

\onecolumn
\section*{Appendix}

\begin{figure*}[h!]
\centering
\includegraphics[width=0.8\textwidth]{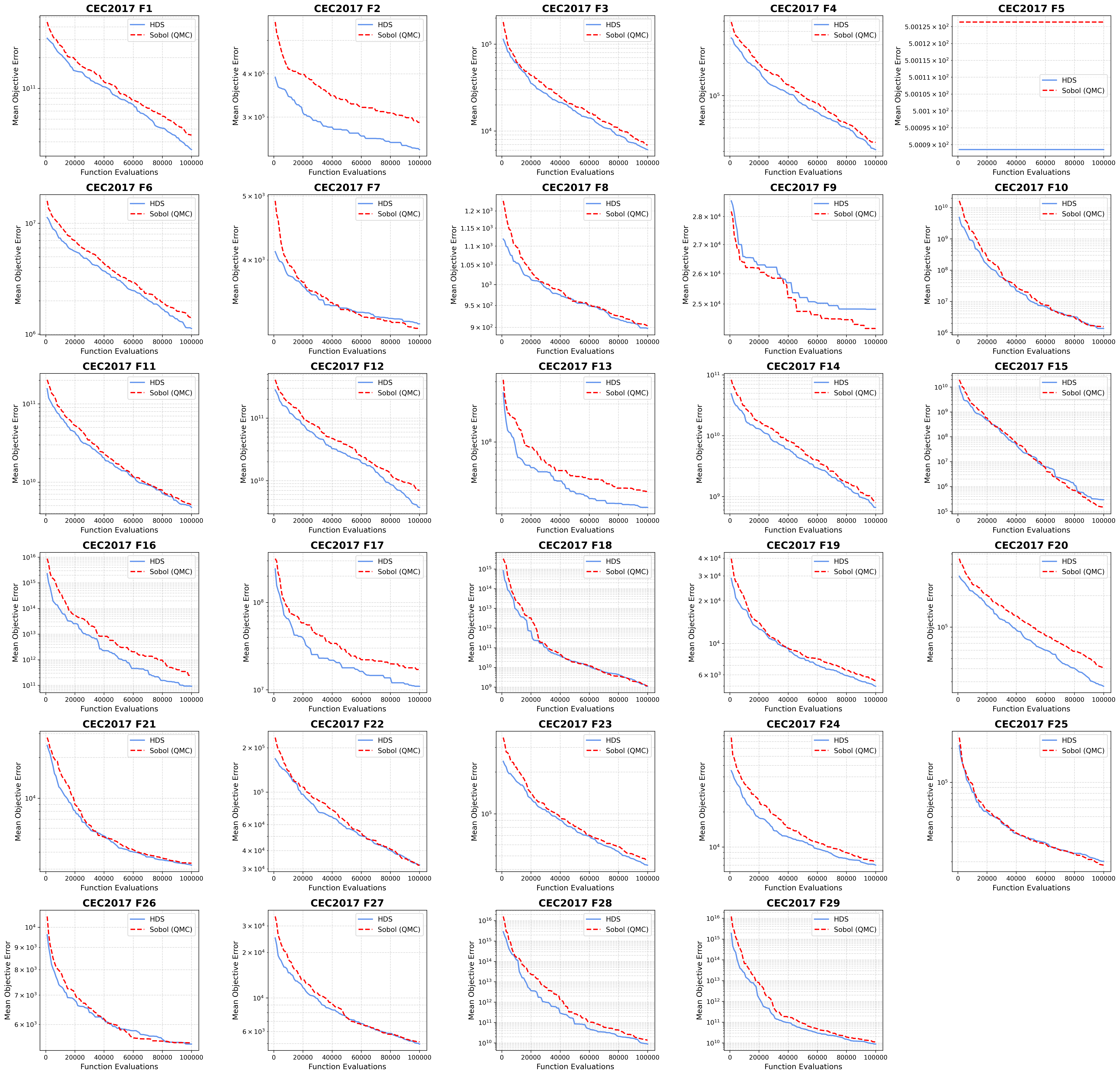}
\caption{Convergence for each test function by number of function evaluations, for HDS (blue) and Sobol (red) ($D=100$, $N=1000$).}
\label{fig:anytime_performance}
\end{figure*}

\newpage
\bibliographystyle{elsarticle-num}
\bibliography{myreferences}

\end{document}